\documentclass[11pt,twoside]{article}

\usepackage{amsmath,amsfonts,amssymb}
\usepackage{graphicx}
\usepackage{eucal}
\usepackage{pxfonts}
\usepackage{mathrsfs}
\usepackage{theorem}
\usepackage{pifont}

\newcommand{\brk}[1]{\left( #1 \right)}
\newcommand{\Brk}[1]{\left[ #1 \right]}
\newcommand{\BRK}[1]{\left\{ #1 \right\}}

\newcommand{\norm}[1]{\left\| #1 \right\|}
\newcommand{\deriv}[2]{\frac{d#1}{d#2}}
\newcommand{\pd}[2]{\frac{\partial#1}{\partial#2}}

\newcommand{\R}{\mathbb{R}}
\newcommand{\calP}{\mathcal{P}}
\newcommand{\calQ}{\mathcal{Q}}
\newcommand{\calI}{\mathcal{I}}
\newcommand{\calM}{\mathcal{M}}
\newcommand{\calC}{\mathcal{C}}

\newcommand{\bs}[1]{\boldsymbol{#1}}
\newcommand{\stre}{\bs{\sigma}}
\newcommand{\vel}{\bs{u}}
\newcommand{\grad}{\bs{\nabla}}
\newcommand{\gvel}{\grad\vel}

\newcommand{\Ltwo}{L^2}
\newcommand{\Linf}{L^{\infty}}
\newcommand{\Linfplus}{L^{\infty,+}}

\newcommand{\Null}{\operatorname{\mathcal{N}}}
\newcommand{\dist}{\operatorname{dist}}
\newcommand{\essinf}{\operatornamewithlimits{ess\ inf}}
\newcommand{\esssup}{\operatornamewithlimits{ess\ sup}}

\newtheorem{theorem}{Theorem}[section]
\newtheorem{proposition}{Proposition}[section]
\newtheorem {corollary}{Corollary}[section]
\newtheorem {lemma}{Lemma}[section]
\newenvironment{proof}{{\flushleft \emph{Proof}:}}{\ding{110}}

\newenvironment{example}{{\flushleft \emph{Example}:}}{}

{\theorembodyfont{\rmfamily}
\newtheorem {assumption}{Assumption}
}

\title{Long-time limit  for a class of quadratic
infinite-dimensional dynamical systems inspired by models of
viscoelastic fluids}

\author{Guy Katriel \footnotemark[1]
\and Raz Kupferman  \footnotemark[1] \and Edriss S. Titi
\footnotemark[2] }
\date{October 11, 2007}
\begin{document}

\maketitle
\renewcommand{\thefootnote}{\fnsymbol{footnote}}
\footnotetext[1]{Institute of Mathematics, The Hebrew University,
Jerusalem 91904, Israel} \footnotetext[2]{Department of Computer
Science and Applied Mathematics, Weizmann Institute of Science,
Rehovot 76100, Israel, Department of Mathematics and Department of
Mechanical and Aerospace Engineering, University of California,
Irvine, CA 92697-3875, USA.}
\renewcommand{\thefootnote}{\arabic{footnote}}


\begin{abstract}
We study a class of quadratic, infinite-dimensional
dynamical systems, inspired by models for viscoelastic fluids. We prove that these equations define a  semi-flow on the cone of positive, essentially bounded functions. As time tends to infinity, the solutions tend to an equilibrium manifold in the $L^2$-norm. Convergence to a particular function on the equilibrium manifold is only proved under additional assumptions. We discuss several possible generalizations.
\end{abstract}


\vskip 0.125in



\section{Introduction}
\label{sec:intro}

This paper is concerned with evolution equations of the form
\begin{equation}
\deriv{y}{t} = y\,\calP(a-y),
\qquad
y(x,0) = y_0(x),
\label{eq:system}
\end{equation}
where $y(\cdot,t)$ is an unknown and $a(\cdot)$ is a given
real-valued function, both defined on a measure space
$(\Omega,\mu)$ with finite mass ($\mu(\Omega)<\infty$). The operator
$\calP$ is an orthogonal projection on the Hilbert space
$\Ltwo = \Ltwo(\Omega,\mu)$. We will use the standard notations
$(\cdot,\cdot)$ for the inner product in $L^2$, and $\|\cdot\|_p$
for the $L^p=L^p(\Omega,\mu)$ norms. We denote by $\Linfplus$ the
cone of strictly positive functions in $\Linf=\Linf(\Omega,\mu)$:
\[
\Linfplus=\{ u\in \Linf : \essinf_{x\in\Omega}{u(x)}>0\}.
\]

Equation \eqref{eq:system} is subject to  the following assumptions:

\begin{assumption}\
\label{ass:1}
\begin{enumerate}
\item[(i)]
The operator $\calP: \Ltwo\rightarrow \Ltwo$ is an
orthogonal projection, satisfying $\calP(\Linf)\subset
\Linf$. Its null space, $\Null(\calP)$, is one-dimensional,  spanned
by an essentially positive function $n\in\Linfplus$, which we take to be normalized,
$(n,n) =1$.
\item[(ii)] The function $a(x)\in \Linf$. Without loss of generality, we can assume that
\[
\calP(a)=a.
\]
\end{enumerate}
\end{assumption}

The system \eqref{eq:system} is a toy model inspired by models of viscoelastic fluids. Specifically, the Maxwell constitutive equation for the conformation tensor is \cite{BAH87}
\begin{equation}
\pd{\stre}{t} + (\vel\cdot\grad)\stre  = (\gvel)^T \stre + \stre (\gvel) +
\frac{1}{\lambda}(\stre - I).
\label{eq:maxwell}
\end{equation}
Here $\stre(x,t)$ is the symmetric positive-definite conformation tensor, $\vel(x,t)$ is the velocity field, $I$ is the unit tensor and $\lambda$ is the elastic relaxation time. For polymers suspended in an incompressible solvent in the creeping flow regime, the velocity field is derived from the conformation tensor through the Stokes equations. The mapping $\stre\mapsto\gvel$, is linear and can be explicitly written by means of a Green function,
\[
\gvel(x) = \int_\Omega G_\Omega(x,y) \cdot \text{div}\stre(y)\,dy,
\]
where $\Omega$ is the domain (which may be bounded on not) and $G_\Omega$ is the corresponding Stokes kernel (i.e., the Green function of the Stokes problem). It can easily be shown that this mapping is, in fact, an orthogonal projection
 (see \cite{FHKK07}). Since the terms in equation \eqref{eq:maxwell} that can potentially lead to finite-time blowup are the stretching terms, it is of interest to omit the advection and the relaxation terms (which by themselves are not norm increasing), and consider systems of the form
\begin{equation}
\pd{\stre}{t}   = [\calP(\bs{a}-\stre)]^T \stre + \stre \calP(\bs{a}-\stre).
\label{eq:fhkk}
\end{equation}
Here $\bs{a}$ is an external force field
(see \cite{FHKK07} for more details).
The system \eqref{eq:system} is a one-dimensional scalar toy model, that mimics the dynamics \eqref{eq:fhkk}.

Equation \eqref{eq:system} can also be viewed as an infinite-dimensional
generalization of a Lotka-Volterra system \cite{Vol31}. In Section~\ref{sec:existence}
we prove that
\eqref{eq:system} defines a global (in time) semi-flow on the cone of
positive functions $\Linfplus$ (Theorem~\ref{th:global} in Section~\ref{sec:existence}).
We then proceed to analyze the
long-time behavior of this system. It is clear that every function $y$
satisfying $\calP(y)=\calP(a)$ is an equilibrium solution of
\eqref{eq:system}, and these are the only equilibria in
$\Linfplus$. Our main theorem
asserts that the equilibrium set
\[
\calM := \BRK{y\in \Linfplus: \calP(y)=a}
\]
is the global attractor for all initial data $y_0\in \Linfplus$ (Theorem~\ref{th:2} in Section~\ref{sec:conv1}). The convergence of
$y(\cdot,t)$ to the manifold $\calM$, as $t\to\infty$, is in the $\Ltwo$-norm.
The theorem does not guarantee uniform convergence, nor does it guarantee
that $y(\cdot,t)$ converges to a specific equilibrium in $\calM$. For this to
happen, additional assumptions are made; various situations are
considered in Section~\ref{sec:conv2}. We conclude this paper with a discussion
about open questions and various generalizations.

\section{Global existence}
\label{sec:existence}

We start by establishing the well-posedness of equation \eqref{eq:system}
under Assumption~\ref{ass:1}.
The first step is to show existence and uniqueness of solutions for short times:

\begin{theorem}[Local-in-time existence and uniqueness]
\label{th:local}
Let $y_0\in\Linf$ be given. Then there exist times $T_1,T_2>0$, depending on $y_0$ only, such that \eqref{eq:system} has a unique solution $y\in C^1((-T_1,T_2),\Linf)$.
\end{theorem}

\begin{proof}
Note first that due to Assumption~\ref{ass:1} the operator $\calP$ has the explicit form
\[
\calP(z) = z - (z,n)n.
\]
It is a bounded linear operator $\Linf\to\Linf$ since
\[
\|\calP(z)\|_\infty \le (1 + \mu(\Omega) \|n\|_\infty^2) \,\|z\|_\infty.
\]
We rewrite \eqref{eq:system} as
\[
\deriv{y}{t} = y a - y \calP(y) \equiv F(y).
\]
The short-time existence and uniqueness of solutions follows from Picard's theorem over Banach spaces, provided that $F$ is a locally Lipschitz continuous mapping $\Linf\to\Linf$. This is indeed the case as $\calP$ is a bounded operator, hence it is locally Lipschitz, and the product of locally Lipschitz functions is again locally Lipschitz.
\end{proof}

We then show that solutions that are initially positive remain so at all times:

\begin{proposition}[Positivity]
\label{prop:positive}
Let $y\in C^1((-T_1,T_2);\Linf)$ be a solution of
\eqref{eq:system}, with initial condition $y_0\in \Linfplus$.
Then $y(x,t)$ remains positive, i.e., $\essinf_{x\in\Omega}y(x,t)>
0$, for all $t\in(-T_1,T_2)$. In other words, the cone $\Linfplus$ is an
invariant set for the dynamics.
\end{proposition}

\begin{proof} The positivity follows readily from the fact that the unique solution of \eqref{eq:system}
solves the integral equation
\begin{equation}{\label{ie}}
y(\cdot,t) = y_0\,\exp\brk{\int_0^t \calP(a - y(\cdot,s))\,ds}.
\end{equation}
\end{proof}

The next step is to show that the solution with initial data in $\Linfplus$, as long as it exists, is bounded, uniformly in time, in $\Linf$, by a constant that only depends on the initial data. The proof relies on the fact that the dynamics \eqref{eq:system} subject to Assumption~\ref{ass:1} preserve the natural order among functions. To simplify notations, we define $\calQ := \calI - \calP$ to be the orthogonal complement of the projection $\calP$, namely, $\calQ y = (n,y)\,n$.

\begin{lemma}
\label{lemma:order}
Let $y\in \Linf$ be a non-negative function, $y\ge 0$. Then,
\[
\essinf_{\Omega} \calQ y(\cdot) \ge 0,
\]
with equality if and only if $y=0$.
\end{lemma}

\begin{proof}
The non-negativity of $y$ and the positivity of $n$ implies that
\[
\essinf_{\Omega} \,\,\calQ y(\cdot) = (n,y) \cdot \essinf_{\Omega} n(\cdot) \ge 0.
\]
Since $\essinf n(x)>0$
equality occurs if and only if $(n,y)=0$, i.e., if and only if $y=0$.
\end{proof}

\begin{proposition}[Comparison principle]
\label{prop:comparison}
Let $y,z\in C^1([0,T);\Linfplus)$ be two solutions of \eqref{eq:system} with initial data $y_0,z_0\in\Linfplus$. If $y_0 \ge z_0$ a.e. in $\Omega$ then
\begin{equation}
y(\cdot,t) \ge z(\cdot,t)
\label{eq:order}
\end{equation}
a.e. in $\Omega$ for all $0\le t< T$.
\end{proposition}

\begin{proof}
Let $t_0$ be the supremum of all values of $t\ge0$ for which the statement holds, i.e.,
$y(x,t) \ge z(x,t)$ a.e. in $\Omega$ for all $0\le t\le t_0$ (it is possible that $t_0=0$). If $t_0=\infty$, there is nothing to prove. If $t_0<\infty$, then by definition
\begin{equation}
y(\cdot,t_0) \ge z(\cdot,t_0).
\label{eq:order0}
\end{equation}
It follows, by Lemma \ref{lemma:order}  that
\[
C := \essinf_{\Omega} \calQ(y(\cdot,t_0)-z(\cdot,t_0)) > 0.
\]

We now define the following sets
\[
\begin{aligned}
\Omega_+ &:= \BRK{x\in\Omega: y(x,t_0) - z(x,t_0) > \tfrac{C}{2}} \\
\Omega_- &:= \BRK{x\in\Omega: y(x,t_0) - z(x,t_0) \le \tfrac{C}{2}}.
\end{aligned}
\]
By the continuity of the mappings $t\mapsto y(\cdot,t)$ and $t\mapsto z(\cdot,t)$ from $[0,T)$ to $\Linf$, there exists a time interval $\delta_1>0$  such that
\begin{equation}
y(x,t) > z(x,t) \qquad \text{ for all $t\in[t_0,t_0+\delta_1)$ for a.e. $x\in \Omega_+$}.
\label{eq:order>}
\end{equation}

We then turn our attention to the set $\Omega_-$, where
\begin{equation}
\begin{split}
\esssup_{\Omega_-} \calP(y(\cdot,t_0)-z(\cdot,t_0)) &\le
\esssup _{\Omega_-} \Brk{y(\cdot,t_0)-z(\cdot,t_0)} - \\
& \essinf _{\Omega_-} \calQ(y(\cdot,t_0)-z(\cdot,t_0))
\le \tfrac{C}{2} - C.
\end{split}
\label{eq:order_P}
\end{equation}

By the differentiability of the mappings $t\mapsto \log y(\cdot,t)$ and $t\mapsto \log z(\cdot,t)$ from $[0,T)$ to $\Linf$, there exists for every $\epsilon>0$ a time interval $\delta_2>0$, such that for all $t\in[t_0,t_0+\delta_2)$,
\[
\begin{gathered}
\norm{\log y(\cdot,t) - \log y(\cdot,t_0) - (t-t_0)\calP(a - y(\cdot,t_0))}_\infty < \epsilon(t-t_0) \\
\norm{\log z(\cdot,t) - \log z(\cdot,t_0) - (t-t_0)\calP(a - z(\cdot,t_0))}_\infty < \epsilon(t-t_0)
\end{gathered}
\]
Thus, for $t\in[t_0,t_0+\delta_2)$,
\[
\begin{split}
\log \frac{y(\cdot,t)}{z(\cdot,t)} &\ge \log \frac{y(\cdot,t_0)}{z(\cdot,t_0)}
- (t-t_0) \calP(y(\cdot,t_0) - z(\cdot,t_0)) - 2\epsilon(t-t_0) \\
&\ge - (t-t_0) \calP(y(\cdot,t_0) - z(\cdot,t_0)) - 2\epsilon(t-t_0),
\end{split}
\]
where the last inequality results from \eqref{eq:order0}. Choosing $\epsilon = C/8$ and using \eqref{eq:order_P} we have that for $t\in[t_0,t_0+\delta_2)$,
\begin{equation}
\inf_{\Omega_-} \log \frac{y(\cdot,t)}{z(\cdot,t)} \ge \frac{C}{4}(t-t_0) \ge 0.
\label{eq:order<}
\end{equation}
Taking $\delta = \min(\delta_1,\delta_2)$ and combining \eqref{eq:order>} and \eqref{eq:order<} we obtain that
\[
y(\cdot,t) \ge z(\cdot,t) \qquad \text{ for all $t\in[t_0,t_0+\delta)$}.
\]
Thus, \eqref{eq:order} holds for all $t\in[t_0,t_0+\delta)$ in contradiction with the definition of $t_0$, which concludes the proof.
\end{proof}

The comparison principle guarantees the boundedness of $y(\cdot,t)$:

\begin{proposition}[Boundedness in $\Linf$]
\label{cor:1}
Let $y\in C^1([0,T);\Linfplus)$ be a solution of \eqref{eq:system} with initial data $y_0$. Then there exists a constant $K>0$, given by \eqref{eq:K} and depending on the initial data, such that
\begin{equation}
\sup_{0\le t<T} y(\cdot,t) \le a + K n(x).
\label{eq:boundednorm}
\end{equation}
\end{proposition}

\begin{proof}
Since $\essinf_\Omega n(x)>0$, then there exists, given $y_0$, a constant $K>0$ such that
\[
z(x) \equiv a(x) + K n(x) \ge y_0(x) \qquad\text{a.e. in $\Omega$}.
\]
Specifically, we can choose
\begin{equation}
K = \esssup_{x\in\Omega} \frac{y_0(x) - a(x)}{n(x)}.
\label{eq:K}
\end{equation}
The function $z$ is an equilibrium solution of \eqref{eq:system}, and by the previous proposition $y(\cdot,t) \le z$ for all $0\le t< T$.
\end{proof}

\begin{theorem}[Global existence]
\label{th:global}
Let $y_0\in\Linfplus$ be given. Then  \eqref{eq:system} has a unique solution $y\in C^1([0,\infty),\Linfplus)$.
\end{theorem}

\begin{proof}
This is a direct consequence of the short-time existence and uniqueness (Theorem~\ref{th:local}) and the bound \eqref{eq:boundednorm} for initial data $y_0\in\Linfplus$.
By the continuation theorem for autonomous ODEs, if $T<\infty$ and $[0,T)$ is the maximal time of existence of the solution $y$, then
\[
\limsup_{t\nearrow T^-} \|y(\cdot,t)\|_\infty = \infty.
\]
Since the norm $ \|y(\cdot,t)\|_\infty$ is continuous in time, this violates the bound \eqref{eq:boundednorm}, hence the maximal existence time is infinite.
\end{proof}

\section{Asymptotic convergence of $y(\cdot,t)$ to $\calM$}
\label{sec:conv1}

Having established the global existence and boundedness of solutions to \eqref{eq:system}, we
proceed to study the long-term behavior of these dynamics.
As in the previous section, it is always assumed that system \eqref{eq:system}
satisfies Assumption~\ref{ass:1}.
The first proposition establishes
the existence of an integral of motion:

\begin{proposition}
\label{prop:integrals}
The functional $\Gamma: \Linfplus\rightarrow \R$
defined by
\[
\Gamma(z) := \int_\Omega n(x) \log z(x)\,d\mu(x), \
\]
is an integral of motion, that is, if $y\in
C^1(\R^+;\Linfplus)$ is a solution of
\eqref{eq:system}, then
\[
\Gamma(y(\cdot,t))=\Gamma(y_0)
\]
for all $t\ge 0$.
\end{proposition}

\begin{proof}
Differentiating we get
\[
\deriv{}{t} \Gamma(y(\cdot,t)) =
 \int_\Omega n(x) \frac{\deriv{}{t} y(x,t)}{y(x,t)}\,d\mu(x) =
(n,\calP(a - y(\cdot,t))) = 0,
\]
where the last equality follows from the symmetry of $\calP$ and the fact that
$n\in\Null(\calP)$.
\end{proof}

The next two propositions reveal the ``dissipative" nature of \eqref{eq:system} through the construction of two Lyapunov functionals. Note that by considering the equilibrium, $\tilde{y}(x)=a(x)+\gamma n(x)$ for
sufficiently large $\gamma$, we have
\[
\essinf_{\Omega} \tilde{y}(\cdot)>0,
\quad\text{ and }\quad
\calP(\tilde{y})=a.
\]

\begin{proposition}
\label{prop:va}
Let $y\in C^1(\R^+;\Linfplus)$ be a solution of
\eqref{eq:system} with $\tilde{y}(x)$ defined as above.
Then the ``entropy" functional
\[
V_a[y(\cdot,t)] := \int_\Omega
\tilde{y}(x)\Brk{\frac{y(x,t)}{\tilde{y}(x)} -
\log\frac{y(x,t)}{\tilde{y}(x)}}\,d\mu(x)
\]
is positive and non-increasing in time.
\end{proposition}

\begin{proof}
The positivity of $V_a$ follows from the fact that $z - \log(z)\ge 1$
for $z>0$, and the positivity of $y(x,t)$ and $\tilde{y}(x)$.
Differentiating along trajectories we get
\begin{eqnarray*}
\deriv{}{t} V_a[y(\cdot,t)] &=& \int_\Omega
\frac{y_t(x,t)}{y(x,t)} \Brk{y(x,t) - \tilde{y}(x)}\, d\mu(x)=
\brk{y(\cdot,t) -\tilde{y}, \calP(a-y(\cdot,t))}\\
&=&\brk{\calP(y(\cdot,t) -\tilde{y}), \calP(a-y(\cdot,t))}= -
\|\calP(a-y(\cdot,t))\|_2^2\le 0,
\end{eqnarray*}
where we have used the fact that $\calP$ is an orthogonal
projection and $\calP(\tilde{y})=a$.
\end{proof}

\begin{proposition}
\label{prop:vb}
Let $y\in C^1(\R^+;\Linfplus)$ be a solution of
\eqref{eq:system}. Then the
``energy" functional
\[
V_b[y(\cdot,t)] := \|\calP(y(\cdot,t)-a)\|_2^2.
\]
is non-increasing in time.
\end{proposition}

\begin{proof}
By explicit differentiation along trajectories we get
\[
\begin{split}
\deriv{}{t} V_b[y(\cdot,t)] &= 2\brk{\calP(y(\cdot,t)-a),\calP(y(\cdot,t) \,\calP(a-y(\cdot,t)))} \\
&= - 2\brk{\calP(y(\cdot,t)-a),y(\cdot,t) \,\calP(y(\cdot,t)-a)} \\
&= -2\|y^{1/2}(\cdot,t) \calP(y(\cdot,t)-a)\|_2^2 \le 0,
\end{split}
\]
where we have used the properties of $\calP$ and the positivity of $y$.
\end{proof}

The identification of the two Lyapunov functionals yields immediately the asymptotic convergence  of $y(\cdot,t)$ to the equilibrium manifold $\calM$.

\begin{theorem}
\label{th:2}
Let $y\in C^1(\R^+;\Linfplus)$ be a
solution of \eqref{eq:system}.
Then
\[
\lim_{t\to\infty} \calP(y(\cdot,t)) = a \qquad\text{ in
$\Ltwo$}.
\]
\end{theorem}

\begin{proof}
We need to prove that
\[
\brk{\dist_{L^2}(y(\cdot,t),\calM)}^2 = \|\calP(y(\cdot,t)-a)\|_2^2 =
V_b(y(\cdot,t))
\]
tends to zero as $t\to\infty$. Since the functionals $V_a,V_b$ are both non-negative, bounded from above (Proposiiton~\ref{cor:1})  and non-increasing in time, both must converge to limits as $t\to\infty$. Since, furthermore,
\[
\deriv{}{t} V_a[y(\cdot,t)] = - V_b[y(\cdot,t)],
\]
the limit of $V_b$ must be zero.
\end{proof}

\begin{example}
\label{ex1}
Assume $\mu(\Omega)=1$ and let $\calP$
be the orthogonal projection in $\Ltwo$
to the space of constants, {\it{i.e.}},
\[
(\calP f)(x) = f(x) - \int_\Omega f(x')\,d\mu(x'),
\]
and $a\in \Linf$ satisfies
\[\int_{\Omega}a(x)d\mu(x)=0.
\]
The system \eqref{eq:system} takes the form
\begin{equation}
\pd{}{t} y(x,t)  = y(x,t) \Big( a(x) +\int_\Omega
y(x',t)\,d\mu(x') -y(x,t) \Big),
\label{eq:example}
\end{equation}
with initial condition $y(\cdot,0)=y_0\in \Linfplus$.
Theorem~\ref{th:global} asserts the existence of a global solution
$y\in C^1(\R^+;\Linfplus)$. By Proposition~\ref{cor:1} there exists a constant $K>0$ such that
\[
\sup_{t\ge 0} y(\cdot,t) \le a  + K.
\]
Finally, by Theorem ~\ref{th:2},
\[
\lim_{t\to\infty}  \brk{y(\cdot,t) - \int_\Omega y(x',t)\,d\mu(x')} = a
\qquad\text{ in $\Ltwo$}.
\]
\end{example}

\section{Asymptotic convergence of $y(\cdot,t)$}
\label{sec:conv2}

We now  question under what conditions does  $y(\cdot,t)$ converge, as $t\to\infty$,  to a specific equilibrium in $\calM$. Note that the $L^2$-convergence of $y(\cdot,t)$  can be decomposed into
\[
\lim_{t\to\infty} y(\cdot,t) =
\lim_{t\to\infty}\calP(y(\cdot,t)) + \lim_{t\to\infty} \calQ(y(\cdot,t)),
\]
where
\[
\calQ(y(\cdot,t)) = (n,y(\cdot,,t))\,n.
\]
We have just proved that the first term on the right-hand side converges to $a$.
It remains to verify under what conditions
\begin{equation}
\beta(t) := ((y(\cdot,t),n)
\label{eq:defbeta}
\end{equation}
converges as $t\to\infty$.

Since, on the one hand, $\calM$ consists of functions of the form $a(x) + \alpha\,n(x)$, for some $\alpha\in\R$, and on the other hand, by Proposition~\ref{prop:integrals} the functional $\Gamma(y(\cdot,t))$ is conserved, the existence of a limiting solution  in $\calM$ requires the following assumption:

\begin{assumption}\
\label{ass:2}
There exists some $y^*\in \calM$ such that
\begin{equation}\label{condition}
\int_\Omega n(x)\,\log y_0(x)\,d\mu(x) = \int_\Omega n(x)\,\log
y^*(x)\,d\mu(x).
\end{equation}
\end{assumption}

Assumption \ref{ass:2} is a restriction on the initial conditions
$y_0$. It assumes the existence of a constant $\alpha$ which solves the
equation
\begin{equation}
\label{equations}
\int_\Omega n(x)\,\log \Brk{ a(x)+ \alpha n(x)}\,d\mu(x)=
\int_\Omega n(x)\,\log y_0(x)\,d\mu(x),
\end{equation}
under the constraint that $\essinf_\Omega [a(x) + \alpha n(x)]>0$.

If we define the set $\calC\subset\R$ by
\begin{equation}
\label{c}
\calC = \BRK{ \xi \in\R
: \essinf_{x\in\Omega} \Brk{ a(x)+\xi n(x)} > 0}
\end{equation}
and $\Phi:\calC\rightarrow \R$ by
\begin{equation}
\label{phi}
\Phi(\xi)= \int_\Omega n(x)\,\log \Brk{a(x)+ \xi n(x)}\,d\mu(x),
\end{equation}
then Assumption \ref{ass:2} is equivalent to the statement
\[
\int_\Omega n(x)\,\log y_0(x)\,d\mu(x) \in \Phi({\calC}).
\]
Note that $\calC$ is in fact an unbounded interval, for $\xi\in\calC$ implies that $\xi_1\in\calC$ for all $\xi_1>\xi$.

The next proposition shows  that such an $\alpha$,
if it exists, is unique.

\begin{proposition}
\label{uniqueness} Given an initial data $y_0\in \Linfplus$,
the function $y^*$ satisfying Assumption \ref{ass:2}, if it
exists, is unique.
\end{proposition}

\begin{proof}
Uniqueness follows at once from the fact that
\[
\deriv{}{\alpha}\Phi(\alpha) = \int_\Omega \frac{n^2(x)}{a(x) + \alpha n(x)}\,d\mu(x) >0
\]
for all $\alpha\in\calC$.
\end{proof}

\begin{example}
\label{ex2}
Consider again the example from the previous section. For concreteness set $\Omega=[0,1]$, with $\mu$ the Lebesgue measure and $a(x) = \sin 2\pi x$. Then, since $n\equiv1$, the equilibria in $\calM$ consist of functions of the form
\[
\sin 2\pi x + \alpha,
\]
where $\alpha>1$, i.e., $\calC = (1,\infty)$. For $\alpha\in\calC$,
\[
\Phi(\alpha) = \int_0^1 \log[\sin2\pi x + \alpha]\,dx > -\log 2.
\]
It follows that Assumption~\ref{ass:2} is satisfied if and only if
\[
\int_0^1 \log y_0(x)\,dx > -\log 2.
\]
\end{example}

The following proposition asserts that the convergence of $y(\cdot,t)$ is guaranteed
if the solution remains bounded away from the boundaries of the cone of positive solutions $\Linfplus$.

\begin{proposition}
\label{be}
If
\begin{equation}
\label{ba}
\liminf_{t\to\infty}\,\,\essinf_{x\in\Omega}y(x,t)>0
\end{equation}
then Assumption~\ref{ass:2} is satisfied. Moreover,
\[
\lim_{t\to\infty} \beta(t) = \alpha,
\]
where $\beta(t)$ is given by \eqref{eq:defbeta} and $\alpha$ is the (unique) solution to  \eqref{equations}. Thus, $y(\cdot,t)\to y^*$ in $\Ltwo$, where $y^* = a + \alpha n$.
\end{proposition}

\begin{proof}
Take any sequence of times $t_m$ that is increasing to infinity.
Since $y(\cdot,t)$ is uniformly bounded in $\Linf$ (Proposition~\ref{cor:1}), then  $\beta(t)$ is bounded, and
there exists a
subsequence $t_{m_k}$ such that $\beta(t_{m_k})$ converges to a
limit $\gamma$, hence
\[
\lim_{k\to\infty} \calQ(y(\cdot,t_{m_k})) =  \gamma
n \qquad\text{ in $\Linf$}.
\]
Theorem~\ref{th:2} implies that
\[
\lim_{k\to\infty} \Brk{y(\cdot,t_{m_k}) - \calQ(y(\cdot,t_{m_k}))
- a} = 0 \qquad\text{ in $\Ltwo$},
\]
from which follows that
\[
\lim_{k\to\infty} y(x,t_{m_k}) = a(x) +  \gamma\,n(x)
\]
in $\Ltwo$, and so
it has a sub-subsequence $y(\cdot,t_{m_{k_j}})$ which converges
a.e. in $\Omega$. Note that \eqref{ba}
implies that a.e. $a(x) + \gamma\,n(x)>0$.
This
implies that
\[\lim_{j\to\infty} n(x) \log y(x,t_{m_{k_j}}) = n(x) \log \Brk{a(x) +
\gamma n(x)}
\]
a.e. Moreover, from  \eqref{ba} and the fact that $y(\cdot,t)$ is uniformly bounded we also have
\[
\sup_{t\geq 0} \| \log y(\cdot,t) \|_{\infty}<\infty.
\]
Using Lebesgue's dominated convergence theorem we conclude
that
\[
\lim_{j\to\infty} \int_{\Omega}n(x) \log y(x,t_{m_{k_j}})d\mu(x) = \int_{\Omega}n(x) \log \Brk{a(x) +
\gamma n(x)}d\mu(x).
\]
By Proposition \ref{prop:integrals} we have for all $m$,
\[
\int_\Omega n(x) \log y(x,t_{m})\,d\mu(x) = \int_\Omega n(x)
\log y_0(x)\,d\mu(x),
\]
therefore
\[
\int_\Omega n(x) \log \Brk{a(x) +  \gamma
n(x)}\,d\mu(x) = \int_\Omega n(x) \log y_0(x)\,d\mu(x).
\]
Thus, Assumption~\ref{ass:2} is satisfied  and it follows, by the uniqueness of $y^*$, hence the uniqueness of $\alpha$ in \eqref{equations}, that $\gamma=\alpha$. We have shown that every sequence $\beta(t_m)$
has a subsequence $\beta(t_{m_{k_j}})$ which converges to
$\alpha$. It follows from an elementary theorem of calculus that $\beta(t)$ tends to $\alpha$ as $t\to\infty$.
This completes the proof.
\end{proof}

Note the immediate corollary:

\begin{corollary}
\label{cor:3}
If Assumption \ref{ass:2} does not hold then
\[
\liminf_{t\to\infty}\,\,\essinf_{x\in\Omega}y(x,t)=0.
\]
\end{corollary}

Condition \eqref{ba} is a sufficient condition for $y(\cdot,t)$ to asymptotically converge to an element of $\calM$. The problem is that it is a property of the solution, and it is not clear \emph{a priori} when does it hold. In the remaining part of this section we establish two situations for which \eqref{ba} holds. In the first case $y_0$ has to be sufficiently large in the following sense:

\begin{proposition}
If there exists a constant $K$ such that
\[
y_0(x) > a(x) + Kn(x) > 0
\qquad\text{a.e. in $\Omega$},
\]
then condition \eqref{ba} holds.
\end{proposition}

\begin{proof}
This is an immediate consequence of the fact that $a+Kn$ is a stationary solution of \eqref{eq:system}, and the comparison principle (Proposition~\ref{prop:comparison}).
\end{proof}

The second situation that can be analyzed is when $a$ and $n$ are simple functions, i.e.,
they have the form
\[
a(x)=\sum_{i=1}^m a_i \chi_{\Omega_i}(x),
\qquad
n(x)=\sum_{i=1}^m n_i \chi_{\Omega_i}(x),
\]
where  $\Omega_1,\dots,\Omega_m$ is a measurable disjoint partition of $\Omega$.

\begin{proposition}
\label{im}
If $a$ and $n$ are simple functions then \eqref{ba} holds.
\end{proposition}

{\bfseries Comment:} The implication of this proposition is that \eqref{ba} holds for any finite-dimensional approximation of \eqref{eq:system}. In particular, the solutions to discrete approximations of \eqref{eq:system} with positive initial data always tend to equilibrium solutions as $t\to\infty$.

\begin{proof}
We first prove the proposition for the particular case in which
$y_0(x)=c>0$ (a constant function).
Note that if $y_0$ and $n$  are simple functions with respect to the partition $(\Omega_i)$, then the right hand side of \eqref{eq:system} is also a simple function,
in which case
$y(x,t)$ is a simple function, constant on each of the sets
$\Omega_i$, for all $t>0$. We denote by $y_i(t)$ the restriction of $y(x,t)$ to the set $\Omega_i$.

Let $M$ be a bound on $|y(x,t)|$ (such a bound is guaranteed to exist by Proposition~\ref{cor:1}). Then for all $t\ge 0$,
\[
\begin{split}
\int_\Omega n(x) \log y(x,t) \,d\mu(x) &=
\sum_{i=1}^m \mu(\Omega_i)\, n_i \log y_i(t)\\
&\le
\log M \,\int_\Omega n(x) \,d\mu(x)+ \brk{\min_{1\leq i \leq m} n_i \mu(\Omega_i)} \log
\brk{\inf_{x\in\Omega}y(x,t)}.
\end{split}
\]
On the other hand, by Proposition \ref{prop:integrals}
\[
\int_\Omega n(x) \log y(x,t) \,d\mu(x)=\int_\Omega n(x) \log y_0(x) \,d\mu(x)
= \log c\,\int_\Omega n(x) \,d\mu(x),
\]
hence
\[
\inf_{\Omega}y(\cdot,t)
\ge  \exp\Brk{\frac{(\log c - \log M)\int_\Omega n(x)\,d\mu(x)}
{\min_{1\leq i \leq m} [n_i \mu(\Omega_i)]}} > 0.
\]
This completes the proof in the case of constant initial conditions.
The general case follows at once from the comparison principle,
as any solution with initial data $y_0\in\Linfplus$ can be
bounded from below by the solution for constant initial data $c = \essinf_\Omega y_0(x)$.
\end{proof}

\section{Discussion}
\label{sec:discussion}

We studied a class of quadratic evolution equations, inspired by models of viscoelastic fluids. Motivated by the physical model, we considered initial data in the cone of positive functions. We showed that the cone of positive $\Linf$ functions is an invariant set, and that solutions in this set exist for all times. As $t\to\infty$ the solutions tend,  in the $\Ltwo$-norm,  to the equilibrium manifold  $\calM$. The convergence of solutions to specific equilibria in $\calM$ could, however, only be proved under additional assumptions.

The following points remain open:
(i) Do solutions always tend to a specific equilibrium if Assumption~\ref{ass:2} is satisfied? We were unable to prove it, nor to find a counter example.
(ii) Do solutions converge, as $t\to\infty$, in situations where Assumption~\ref{ass:2} does not hold? While, in such case, the solution cannot converge to an equilibrium in $\calM$ (Corollary~\ref{cor:3}), it can, in principle, converge to an equilibrium on the boundary of the cone,
\[
\overline{\Linfplus} = \BRK{y\in\Linf:\,\, y(x)\ge 0}.
\]
(iii) Does the solution converge to $\calM$ in any $L^p$-norm, for $p>2$, and in particular, for $p=\infty$?

Another question is whether our results remain valid when the kernel of the projection $\calP$ has dimension greater than one. The comparison principle (Proposition~\ref{prop:comparison}) no longer holds in this case, and as a result, we no longer have a bound on the $\Linf$ norm, nor do we have a global existence theorem. Assuming, however, that a solution does exist for all times, it is easy to see that Proposition~\ref{prop:vb} still holds, i.e., the ``energy" functional $V_b$ is a Lyapunov functional. To prove that the ``entropy" functional $V_a$ is also a Lyapunov functional, we need to have a positive function $\tilde{y}$ such that $\calP(\tilde{y}) = a$. If such function exists then Proposition~\ref{prop:va} remains valid, and $\calP(y)$ tends to $a$ in the $\Ltwo$-norm (Theorem~\ref{th:2}).

System \eqref{eq:system} can be generalized in many different ways, for example, with $y$ being a matrix valued function and products reinterpreted as matrix products; this is indeed the appropriate setting in the viscoelastic context \cite{FHKK07}. Another generalization of \eqref{eq:system} is when $\calP$ is a general non-negative operator (not necessarily a projection), i.e., $(y,\calP(y)) \ge 0$ for all $y\in\Ltwo$. We believe that such a system still exhibits global-in-time existence for positive initial data, as well as asymptotic convergence.

\bigskip\noindent
{\bfseries Acknowledgments}
We are grateful to Raanan Fattal for discussions that motivated this present work.
GK was partially supported
by the Edmund Landau
Center for Research in Mathematical Analysis and Related Areas,
sponsored by the Minerva Foundation (Germany).
RK was partially supported by the Israel Science
Foundation founded by the Israel Academy of Sciences and Humanities,
and by the Applied Mathematical
Sciences subprogram of the Office of Energy Research of the US
Department of Energy under Contract DE-AC03-76-SF00098.
The work of EST was supported in part by
the NSF grant no.~DMS-0504619, the ISF grant no.~120/6, and the BSF
grant no.~2004271.


\begin{thebibliography}{1}

\bibitem{BAH87}
{\sc R.~Bird, R.~Armstrong, and O.~Hassager}, {\em Dynamics of polymeric
  liquids. {Volume 1}}, John Wiley and Sons, New York, 1987.

\bibitem{FHKK07}
{\sc R.~Fattal, O.~Hald, G.~Katriel, and R.~Kupferman}, {\em Global stability
  of equilibrium manifolds, and "peaking" behavior in quadratic differential
  systems related to viscoelastic models}, J. Non-Newton. Fluid Mech., 144
  (2007), pp.~30--41.

\bibitem{Vol31}
{\sc V.~Volterra}, {\em Variations and fluctuations of the number of
  individuals in animal species living together}, in Animal ecolocy,
  McGraw-Hill, 1931.

\end{thebibliography}
\end{document}